\magnification=1200
%\nopagenumbers

\hsize=11.25cm    
\vsize=18cm       
\parindent=12pt   \parskip=5pt     

\hoffset=.5cm   
\voffset=.8cm   

\pretolerance=500 \tolerance=1000  \brokenpenalty=5000

\catcode`\@=11

\font\eightrm=cmr8         \font\eighti=cmmi8
\font\eightsy=cmsy8        \font\eightbf=cmbx8
\font\eighttt=cmtt8        \font\eightit=cmti8
\font\eightsl=cmsl8        \font\sixrm=cmr6
\font\sixi=cmmi6           \font\sixsy=cmsy6
\font\sixbf=cmbx6

\font\tengoth=eufm10 
\font\eightgoth=eufm8  
\font\sevengoth=eufm7      
\font\sixgoth=eufm6        \font\fivegoth=eufm5

\skewchar\eighti='177 \skewchar\sixi='177
\skewchar\eightsy='60 \skewchar\sixsy='60

\newfam\gothfam           \newfam\bboardfam

\def\tenpoint{
  \textfont0=\tenrm \scriptfont0=\sevenrm \scriptscriptfont0=\fiverm
  \def\rm{\fam\z@\tenrm}
  \textfont1=\teni  \scriptfont1=\seveni  \scriptscriptfont1=\fivei
  \def\oldstyle{\fam\@ne\teni}\let\old=\oldstyle
  \textfont2=\tensy \scriptfont2=\sevensy \scriptscriptfont2=\fivesy
  \textfont\gothfam=\tengoth \scriptfont\gothfam=\sevengoth
  \scriptscriptfont\gothfam=\fivegoth
  \def\goth{\fam\gothfam\tengoth}
  
  \textfont\itfam=\tenit
  \def\it{\fam\itfam\tenit}
  \textfont\slfam=\tensl
  \def\sl{\fam\slfam\tensl}
  \textfont\bffam=\tenbf \scriptfont\bffam=\sevenbf
  \scriptscriptfont\bffam=\fivebf
  \def\bf{\fam\bffam\tenbf}
  \textfont\ttfam=\tentt
  \def\tt{\fam\ttfam\tentt}
  \abovedisplayskip=12pt plus 3pt minus 9pt
  \belowdisplayskip=\abovedisplayskip
  \abovedisplayshortskip=0pt plus 3pt
  \belowdisplayshortskip=4pt plus 3pt 
  \smallskipamount=3pt plus 1pt minus 1pt
  \medskipamount=6pt plus 2pt minus 2pt
  \bigskipamount=12pt plus 4pt minus 4pt
  \normalbaselineskip=12pt
  \setbox\strutbox=\hbox{\vrule height8.5pt depth3.5pt width0pt}
  \let\bigf@nt=\tenrm       \let\smallf@nt=\sevenrm
  \normalbaselines\rm}

\def\eightpoint{
  \textfont0=\eightrm \scriptfont0=\sixrm \scriptscriptfont0=\fiverm
  \def\rm{\fam\z@\eightrm}
  \textfont1=\eighti  \scriptfont1=\sixi  \scriptscriptfont1=\fivei
  \def\oldstyle{\fam\@ne\eighti}\let\old=\oldstyle
  \textfont2=\eightsy \scriptfont2=\sixsy \scriptscriptfont2=\fivesy
  \textfont\gothfam=\eightgoth \scriptfont\gothfam=\sixgoth
  \scriptscriptfont\gothfam=\fivegoth
  \def\goth{\fam\gothfam\eightgoth}
  
  \textfont\itfam=\eightit
  \def\it{\fam\itfam\eightit}
  \textfont\slfam=\eightsl
  \def\sl{\fam\slfam\eightsl}
  \textfont\bffam=\eightbf \scriptfont\bffam=\sixbf
  \scriptscriptfont\bffam=\fivebf
  \def\bf{\fam\bffam\eightbf}
  \textfont\ttfam=\eighttt
  \def\tt{\fam\ttfam\eighttt}
  \abovedisplayskip=9pt plus 3pt minus 9pt
  \belowdisplayskip=\abovedisplayskip
  \abovedisplayshortskip=0pt plus 3pt
  \belowdisplayshortskip=3pt plus 3pt 
  \smallskipamount=2pt plus 1pt minus 1pt
  \medskipamount=4pt plus 2pt minus 1pt
  \bigskipamount=9pt plus 3pt minus 3pt
  \normalbaselineskip=9pt
  \setbox\strutbox=\hbox{\vrule height7pt depth2pt width0pt}
  \let\bigf@nt=\eightrm     \let\smallf@nt=\sixrm
  \normalbaselines\rm}

\tenpoint

\def\pc#1{\bigf@nt#1\smallf@nt}         \def\pd#1 {{\pc#1} }

\catcode`\;=\active
\def;{\relax\ifhmode\ifdim\lastskip>\z@\unskip\fi
\kern\fontdimen2  -1.2 \fontdimen3 \string;}

\catcode`\:=\active
\def:{\relax\ifhmode\ifdim\lastskip>\z@\unskip\fi\penalty\@M\ \fi\string:}

\catcode`\!=\active
\def!{\relax\ifhmode\ifdim\lastskip>\z@
\unskip\fi\kern\fontdimen2  -1.1 \fontdimen3 \string!}

\catcode`\?=\active
\def?{\relax\ifhmode\ifdim\lastskip>\z@
\unskip\fi\kern\fontdimen2  -1.1 \fontdimen3 \string?}

\frenchspacing

\def\raggedbottom{\topskip 10pt plus 36pt\r@ggedbottomtrue}

\def\pointir{\unskip . --- \ignorespaces}

\def\Medbreak{\vskip-\lastskip\medbreak}

\long\def\th#1 #2\enonce#3\endth{
   \Medbreak\noindent
   {\pc#1} {#2\unskip}\pointir{\it #3}\smallskip}

\def\decale#1{\smallbreak\hskip 28pt\llap{#1}\kern 5pt}
\def\decaledecale#1{\smallbreak\hskip 34pt\llap{#1}\kern 5pt}
\def\puce{\smallbreak\hskip 6pt{$\scriptstyle\bullet$}\kern 5pt}

\def\eqalign#1{\null\,\vcenter{\openup\jot\m@th\ialign{
\strut\hfil$\displaystyle{##}$&$\displaystyle{{}##}$\hfil
&&\quad\strut\hfil$\displaystyle{##}$&$\displaystyle{{}##}$\hfil
\crcr#1\crcr}}\,}

\catcode`\@=12

\showboxbreadth=-1  \showboxdepth=-1

\newcount\numerodesection \numerodesection=1
\def\section#1{\bigbreak
 {\bf\number\numerodesection.\ \ #1}\nobreak\medskip
 \advance\numerodesection by1}

\mathcode`A="7041 \mathcode`B="7042 \mathcode`C="7043 \mathcode`D="7044
\mathcode`E="7045 \mathcode`F="7046 \mathcode`G="7047 \mathcode`H="7048
\mathcode`I="7049 \mathcode`J="704A \mathcode`K="704B \mathcode`L="704C
\mathcode`M="704D \mathcode`N="704E \mathcode`O="704F \mathcode`P="7050
\mathcode`Q="7051 \mathcode`R="7052 \mathcode`S="7053 \mathcode`T="7054
\mathcode`U="7055 \mathcode`V="7056 \mathcode`W="7057 \mathcode`X="7058
\mathcode`Y="7059 \mathcode`Z="705A

% handling accented characters in plain TeX :

\def\hfl#1#2#3{\smash{\mathop{\hbox to#3{\rightarrowfill}}\limits
^{\textstyle#1}_{\textstyle#2}}}

\def\Q{{\bf Q}}

\def\R{{\bf R}}
\def\C{{\bf C}}

\def\Z{{\bf Z}}

\def\F{{\bf F}}
\def\Fp{{\F_{\!p}}}

\def\Aut{\mathop{\rm Aut}\nolimits}
\def\GL{\mathop{\rm GL}\nolimits}

\def\Id{\mathop{\rm Id}\nolimits}

\def\Gal{\mathop{\rm Gal}\nolimits}
\def\Frob{\mathop{\rm Frob}\nolimits}

\def\to{\rightarrow}

\def\normressym(#1,#2)_#3{\displaystyle\left({#1,#2\over#3}\right)}

\def\mod{\mathop{\rm mod.}\nolimits}
\def\pmod#1{\;(\mod#1)}

\newcount\refno 
\long\def\ref#1:#2<#3>{                                        
\global\advance\refno by1\par\noindent                              
\llap{[{\bf\number\refno}]\ }{#1} \pointir{\it #2} #3\goodbreak }

\def\citer#1(#2){[{\bf\number#1}\if#2\empty\relax\else,\ {#2}\fi]}

\newbox\bibbox
\setbox\bibbox\vbox{
\bigskip
\centerline{---$*$---$*$---}
\bigbreak
\centerline{{\pc REFERENCES}}

\ref{\pc ARTHUR} (J):
Automorphic representations and number theory,
<in CMS Conf.\ Proc., 1, Amer.\ Math.\ Soc., Providence, 1981.>
\newcount\arthur \global\arthur=\refno

\ref{\pc ARTHUR} (J):
The principle of functoriality,
<Bull.\ Amer.\ Math.\ Soc.\ (N.S.) {\bf 40} (2003) 1, 39--53.>
\newcount\arthurmmiii \global\arthurmmiii=\refno

\ref{\pc CASSELS} (J) \& {\pc FR{\"O}HLICH} (A) (Eds.):
Algebraic number theory.
<Proceedings of the instructional conference held at the University of Sussex,
Brighton, September 1--17, 1965, Corrected reprint, London Mathematical
Society, 2010.>
\newcount\ant \global\ant=\refno

\ref{\pc COX} (D):
Primes of the form $x^2 + ny^2$. Fermat, class field theory and complex
multiplication.  
<John Wiley \& Sons, New York, 1989.>
\newcount\cox \global\cox=\refno

\ref{\pc GELBART} (S):
An elementary introduction to the Langlands program,
<Bull.\ Amer.\ Math.\ Soc.\ (N.S.) {\bf 10}  (1984) 2, 177--219.>
\newcount\gelbart \global\gelbart=\refno

\ref{\pc HARDER} (G):
The Langlands program (an overview),
<School on Automorphic Forms on ${\rm GL}(n)$, 207--235,
ICTP Lect.\ Notes, 21, Abdus Salam Int.\ Cent.\ Theoret.\ Phys., Trieste, 2008.>
\newcount\harder \global\harder=\refno

\ref{\pc HONDA} (T):
Invariant differentials and $L$-functions. Reciprocity law for quadratic
fields and elliptic curves over $\Q$, 
<Rend.\ Sem.\ Mat.\ Univ.\ Padova {\bf 49} (1973), 323--335.>
\newcount\honda \global\honda=\refno

\ref{\pc KNAPP} (A):
Introduction to the Langlands program,
<in Proc.\ Sympos.\ Pure Math., 61, Amer.\ Math.\ Soc., Providence, 1997.>
\newcount\knapp \global\knapp=\refno

\ref{\pc LANGLANDS} (R):
Representation theory: its rise and its role in number theory,
<Proceedings of the Gibbs Symposium (New Haven, 1989), 181--210,
Amer.\ Math.\ Soc., Providence, RI, 1990.>
\newcount\gibbs \global\gibbs=\refno

\ref{\pc SERRE} (J-P):
Modular forms of weight one and Galois representations,
<in  Proc. Sympos. (Durham, 1975),  193--268, Academic Press, London, 1977.> 
\newcount\serremodi \global\serremodi=\refno

\ref{\pc SERRE} (J-P):
On a theorem of Jordan,
<Bull.\ Amer.\ Math.\ Soc.\ (N.S.) {\bf 40} (2003), no.\ 4, 429--440.>
\newcount\serre \global\serre=\refno

\ref{\pc SHIMURA} (G):
A reciprocity law in non-solvable extensions,
<J.\ Reine Angew.\ Math.\ {\bf 221} (1966) 209--220.>
\newcount\shimura \global\shimura=\refno

\ref{\pc TAYLOR} (R):
 Reciprocity laws and density theorems,
<on the author's website.>
\newcount\taylor \global\taylor=\refno

\ref{\pc WYMAN} (B):
What is a reciprocity law~?
<Amer.\ Math.\ Monthly {\bf 79} (1972), 571--586; {\it correction}, ibid.\
{\bf 80} (1973), 281.>
\newcount\wyman \global\wyman=\refno

\ref{\pc ZAGIER} (D):
Modular points, modular curves, modular surfaces and modular forms,
<Lecture Notes in Math., 1111, Springer, Berlin, 1985.>
\newcount\zagier \global\zagier=\refno

}

\centerline{\bf Splitting primes}

\vskip1cm
\centerline{Chandan Singh Dalawat}
\medskip

{\eightpoint We give an elementary introduction, through illustrative examples
  but without proofs, to one of the basic consequences of the Langlands
  programme, namely the law governing the primes modulo which a given
  irreducible integral polynomial splits completely.  Some recent results,
  such as the modularity of elliptic curves over the rationals, or the proof
  of Serre's conjecture by Khare and Wintenberger, are also illustrated
  through examples.}

\bigskip

{\it Certainly the best times were when I was alone with mathematics, free of
  ambition and pretense, and indifferent to the world.} --- Robert Langlands.

\bigskip

Groups, rings, fields, polynomials, primes, functions~: most mathematics
students at the university have come across these concepts.  It so happens
that one of the basic problems of number theory can be formulated using
nothing more than these notions~; we shall try to explain Langlands'
conjectural solution to this problem.  Everything we say here has been said
before, but some readers might not have access to the original sources.  This
is the pretext for putting together these observations~; the aim is not so
much to instruct as to convey the impression that these examples are part of a
larger pattern. Proofs are altogether omitted.

\smallbreak

Let us consider a monic polynomial $f=T^n+c_{n-1}T^{n-1}+\cdots+c_1T+c_0$ of
degree $n>0$ with coefficients $c_i$ in the ring $\Z$ of rational integers
(the adjective {\it monic} means that the coefficient of $T^n$ is~$1$).
Suppose that the polynomial $f$ is irreducible~: it cannot be written as a
product $f=gh$ of two polynomials $g,h\in\Z[T]$ of degree $<n$.  Basic
examples to be kept in mind are $f=T^2+1$ or $f=T^3-T-1$.  (There are
irreducible polynomials of every possible degree.  For example, Selmer showed
that $T^n-T-1$ is irreducible for every $n>1$.)

For every number $m>0$, we have the finite ring $\Z/m\Z$ with $m$ elements
each of which is of the form $\bar a$ for some $a\in\Z$, with the relation
$\bar a=\bar b$ if $m$ divides $a-b$~; this relation is also written $a\equiv
b\pmod m$, and we say that $\bar a$ is the {\it reduction\/} of $a$
modulo~$m$.  The group of invertible elements of this ring is denoted
$(\Z/m\Z)^\times$, so that $\bar a\in(\Z/m\Z)^\times$ if and only if
$\gcd(a,m)=1$.  The order of this group is denoted $\varphi(m)$.  When the
integer $m$ is a prime number~$p$, the ring $\Z/p\Z$ is a field denoted $\Fp$
and the group $\F_p^\times$ is cyclic of order $\varphi(p)=p-1$.

For every prime~$p$, the polynomial $f\in\Z[T]$ gives rise to a degree-$n$
polynomial $\bar f=T^n+\bar c_{n-1}T^{n-1}+\cdots+\bar c_1T+\bar c_0$ whose
coefficients $\bar c_i\in\Fp$ are the reductions of $c_i$ modulo~$p$~; this
$\bar f\in\Fp[T]$ is called the {\it reduction\/} of $f$ modulo~$p$.

Now, although our polynomial $f$ is irreducible by hypothesis, there may be
some primes $p$ modulo which its reduction $\bar f$ has $n$ distinct roots.
For example, if $f=T^2+1$ and $p=5$, then $\bar f=(T+\bar 2)(T-\bar 2)$.  More
generally, we have $\bar f=(T+\bar1)^2$ if $p=2$, $\bar f$ has two distinct
roots if $p\equiv1\pmod4$ and $\bar f$ remains irreducible if
$p\equiv-1\pmod4$.  One sees this by remarking that, for odd~$p$, a root
$x\in\Fp$ of $\bar f$ is an order-$4$ element in the cyclic group
$\F_p^\times$ of order $p-1$, so that $x$ exists if and only if $4\,|\,p-1$.
This example goes back to Fermat.

Take $f=T^2-T-1$ as the next example.  If $p=5$, then $\bar f=(T+\bar2)^2$.
It can be checked that for $p\neq5$, the polynomial $\bar f$ has two distinct
roots if and only if $p\equiv\pm1\pmod5$, and that it is irreducible if and
only if $p\equiv\pm2\pmod 5$.

We notice that for these two quadratic polynomials $f$, the primes $p$ for
which $\bar f$ has two distinct roots are specified by ``congruence
conditions''~: by $p\equiv1\pmod4$ when $f=T^2+1$, by $p\equiv\pm1\pmod5$ when
$f=T^2-T-1$.

{\it It can be shown that for any (monic, irreducible) quadratic polynomial\/
  $f\in\Z[T]$, the primes\/ $p$ for which\/ $\bar f\in\Fp[T]$ has two distinct
  roots are given by certain congruence conditions modulo some number\/ $D_f$
  depending on\/ $f$.}  This statement (with a precise expression for $D_f$
which we omit) is implied by the law of quadratic reciprocity, first proved by
Gauss.  The law in question says that for any two odd primes $p\neq q$, we
have
$$
\left(p\over q\right)\left(q\over p\right)=(-1)^{{p-1\over2}{q-1\over2}}
\leqno{(1)}
$$
where, by definition, $\left(p\over q\right)=1$ if $\bar p\in\F_q^\times$
is a square (if $\bar p=\bar a^2$ for some $\bar a\in\F_q^\times$), and
$\left(p\over q\right)=-1$ if $\bar p\in\F_q^\times$ is not a square.

With the help of this law, the reader should be able to check that the
reduction $\bar f$ (modulo $p$) of $f=T^2-11$ has two distinct roots if and
only if, modulo~$44$,
$$
\bar p\in\{\bar1,\bar5,\bar7,\bar9,
\overline{19},\overline{25},\overline{35},\overline{37},
\overline{39},\overline{43}\}.
$$
(Incidentally, given any integer $m>0$ and any $\bar a\in(\Z/m\Z)^\times$,
there are infinitely many primes $p\equiv a\pmod m$~; in a precise sense, the
proportion of such primes is $1/\varphi(m)$, as shown by Dirichlet.)

The law of quadratic reciprocity can be formulated as an equality of two
different kinds of $L$-{\it functions}, an ``Artin $L$-function'' and a
``Hecke $L$-function'', both functions of a complex variable $s$.  The first
one carries information about the primes $p$ modulo which the polynomial
$T^2-q^*$ (where $q^*=(-1)^{(q-1)/2}q$, so that $q^*\equiv1\pmod4$) has two
distinct roots, and is defined as
$$
L_1(s)=\prod_p{1\over 1-\left({q^*\over p}\right)p^{-s}}
$$
where $\left({q^*\over 2}\right)=(-1)^{(q^2-1)/8}$ and $\left({q^*\over
    q}\right)=0$. The second one carries information about  which primes~$p$
lie in which arithmetic progressions modulo~$q$, and is defined as
$$
L_2(s)%=\sum_{n>0}\left({n\over q}\right)n^{-s}
=\prod_p{1\over 1-\left({p\over q}\right)p^{-s}},
$$
where we put $\left({q\over q}\right)=0$.  The law of quadratic reciprocity
$(1)$ is equivalent to the statement that $L_1(s)=L_2(s)$ (except possibly for
the factor at the prime~$2$), as can be seen by comparing the coefficients of
$p^{-s}$ and noting that $\left({-1\over p}\right)=(-1)^{(p-1)/2}$.  

\medbreak

What are the (monic, irreducible) degree-$n$ polynomials $f$ for which the
primes $p$ where the reduction $\bar f$ has~$n$ distinct roots in $\Fp$ can be
specified by similar congruence conditions~?  Before characterising such
polynomials, let us discuss one more example.

For every prime $l$, consider the (unique) monic polynomial $\Phi_l\in\C[T]$
whose roots are precisely the $l-1$ primitive $l$-th roots $e^{2i\pi k/l}$
($0<k<l$) of~$1$, namely $\Phi_l=T^{l-1}+T^{l-2}+\cdots+T+1$~; notice that
$\Phi_l\in\Z[T]$.  It can be shown (in much the same way as our discussion of
$T^2+1$ above) that the reduction $\bar\Phi_l$ modulo some prime $p\neq l$ has
$l-1$ distinct roots in $\Fp$ if and only if $l\,|\,p-1$, or equivalently
$p\equiv1\pmod l$.

More generally, for every integer $m>0$, the primitive $m$-th roots of~$1$ in
$\C$ are $\zeta_m^k$, where $\zeta_m=e^{2i\pi/m}$ and $\gcd(k,m)=1$.  If
$k'\equiv k\pmod m$, then $\zeta_m^{k'}=\zeta_m^k$, so the number of such
primitive roots is $\varphi(m)$, the order of the group $(\Z/m\Z)^\times$ (as
$\gcd(k,m)=1$ if and only if $\bar k\in(\Z/m\Z)^\times$).  Let
$\Phi_m\in\C[T]$ be the monic polynomial whose roots are precisely these
$\varphi(m)$ numbers $\zeta_m^k$, where $\gcd(k,m)=1$, so that $\Phi_4=T^2+1$.
It can be shown that $\Phi_m\in\Z[T]$ and that it is irreducible~; $\Phi_m$ is
called {\it the cyclotomic polynomial\/} of level~$m$.

For which primes $p$ does the reduction $\bar\Phi_m$ have $\varphi(m)$
distinct roots in $\Fp$~?  Some work is required to answer this question, and
it must have been first done by Kummer in the XIX$^{\rm th}$ century, if not
by Gauss before him.  {\it It turns out that\/ the reduction\/ $\bar\Phi_m$
  modulo\/~$p$ of the level-$m$ cyclotomic polynomial\/ $\Phi_m$ has\/
  $\varphi(m)=\deg\Phi_m$ distinct roots in\/ $\Fp$ if and only if\/
  $p\equiv1\pmod m$.}

Now, for any monic irreducible polynomial $f$ of degree~$n>0$, we can adjoin
the $n$ roots of $f$ in $\C$ to $\Q$ to obtain a field $K_f$.  The dimension
of $K_f$ as a vector space over $\Q$ is finite (and divides $n!$).  For
example, when $f=T^2+1$ and $i\in\C$ is a root of $f$, the field $K_f$
consists of all complex numbers $x+iy$ such that $x,y\in\Q$.  When $f=\Phi_l$
for some prime~$l$, $K_f$ consists of all $z\in\C$ which can be written
$z=b_0+b_1\zeta_l+\cdots+b_{l-2}\zeta_l^{l-2}$ for some $b_i\in\Q$.

The automorphisms of the field $K_f$ form a group which we denote by $\Gal_f$
in honour of Galois.  When $f$ is quadratic, $\Gal_f$ is isomorphic to
$\Z/2\Z$~; when $f=\Phi_m$ is the cyclotomic polynomial of level $m>0$, then
$\Gal_f$ is canonically isomorphic to $(\Z/m\Z)^\times$.  The group $\Gal_f$
is commutative in these two cases~; we say that the polynomial $f$ is {\it
  abelian\/} (in honour of Abel) if the group $\Gal_f$ is commutative.

{\it It has been shown that, $f\in\Z[T]$ being an irreducible monic polynomial
  of degree~$n>0$, the set of primes~$p$ modulo which\/ $\bar f$ has\/ $n$
  distinct roots can be characterised by congruence conditions modulo some
  number\/ $D_f$ depending only on\/ $f$ if and only if the group\/ $\Gal_f$
  is commutative (as is for example the case when $n=2$ or\/ $f=\Phi_m$ for\/
  some $m>0$).}

Visit {\tt http\string:/$\!$/mathoverflow.net/questions/11688} for a
discussion of this equivalence~; the precise version of this statement is a
consequence of the following beautiful theorem of Kronecker and Weber~: {\it
  For every abelian polynomial\/ $f$, there is some integer\/ $m>0$ such
  that\/ $K_f\subset K_{\Phi_m}$~: every root of\/ $f$ can be written as a
  linear combination of\/ $1,\zeta_m,\ldots,\zeta_m^{\varphi(m)-1}$ with
  coefficients in $\Q$.}

The generalisation of this theorem to all ``global fields'' constitutes
``class field theory'', which was developed in the first half of the XX$^{\rm
  th}$ century by Hilbert, Furtw{\"a}ngler, Takagi, Artin, Chevalley, Tate,
$\ldots$ This theory is to abelian polynomials $f$ what the quadratic
reciprocity law is to quadratic polynomials over $\bf Q$.  As we saw above for
the latter, class field theory can be formulated as the equality between an
``Artin $L$-function'' and a ``Hecke $L$-function'', the first one carrying
information about the primes modulo which $f$ splits completely and the second
carrying information about which primes lie in which arithmetic progressions
modulo~$D_f$.

How about polynomials $f$ which are not abelian, the ones for which the group
$\Gal_f$ is not commutative~?  Such an $f$, namely $f=T^3-2$, was already
considered by Gau\ss.  He proved that this $f$ splits completely modulo a
prime $p\neq2,3$ if and only if $p=x^2+27y^2$ for some $x,y\in\Z$.  For
general~$f$, how can we characterise the primes $p$ modulo which $\bar f$ has
$n=\deg f$ distinct roots~?  The Langlands programme, initiated in the late
60s, provides the only satisfactory, albeit conjectural, answer to this basic
question (in addition to many other basic questions).  This is what we would
like to explain next.

Denote by $N_p(f)$ the number of roots of $\bar f$ in $\Fp$, so that for the
abelian polynomial $f=T^2-T-1$ we have $N_5(f)=1$ and
$$
N_p(f)=\cases{2&if $p\equiv\pm1\pmod5$,\cr
0&if $p\equiv\pm2\pmod5$,\cr}
$$
as we saw above.  It can be checked that $N_p(f)=1+a_p$ for all
primes $p$, where $a_n$ is the coefficient of $q^n$ in the formal power series
$$
{q-q^2-q^3+q^4\over 1-q^5}=q-q^2-q^3+q^4+q^6-q^7-q^8+q^9+\cdots
$$
in the indeterminate $q$.  We can be said to have found a ``formula'' for
$N_p(T^2-T-1)$. Notice that the $a_n$ are {\it strongly multiplicative\/} in
the sense that $a_{mm'}=a_ma_{m'}$ for all $m>0$, $m'>0$.

It follows from class field theory that there are similar ``formulae'' or
``reciprocity laws'' for $N_p(f)$ for every abelian polynomial $f$.  What
about polynomials $f$ for which $\Gal_f$ is {\it not\/} abelian~?  

The fundamental insight of Langlands was that there is a ``formula'' for
$N_p(f)$ for every polynomial $f$, abelian or not, and to say precisely what
is meant by ``formula''.  We cannot enter into the details but content
ourselves with providing some illustrative examples (without proof).  I learnt
the first one, which perhaps goes back to Hecke, from Emerton's answer to one
of my questions (visit {\tt
  http\string:/$\!$/mathoverflow.net/questions/12382}) and also from a paper
of Serre \citer\serre().

\medbreak

Consider the polynomial $f=T^3-T-1$ which is {\it not\/} abelian~; the group
$\Gal_f$ is isomorphic to the smallest group which is not commutative, namely
the symmetric group ${\goth S}_3$ of all permutations of three letters.
Modulo $23$, $\bar f$ has a double root and a simple root, so $N_{23}(f)=2$.
If a prime $p\neq23$ is a square in $\F_{23}^\times$, then $p$ can be written
either as $x^2+xy+6y^2$ or as $2x^2+xy+3y^2$ (but not both) for some
$x,y\in\Z$~; we have $N_p(f)=3$ in the first case whereas $N_p(f)=0$ in the
second case.  Finally, if a prime $p\neq23$ is not a square in
$\F_{23}^\times$, then $N_p(f)=1$.  The smallest $p$ for which $N_p(f)=3$ is
$59=5^2+5.2+6.2^2$. 

It can be checked that for all primes $p$, we have $N_p(f)=1+a_p$, where
$a_n$ is the coefficient of $q^n$ in the formal power series
$\eta_{1^1,23^1}$~:  
$$
q\prod_{k=1}^{+\infty}(1-q^k)(1-q^{23k})
=q-q^2-q^3+q^6+q^8-q^{13}-q^{16}+q^{23}-q^{24}+\cdots
$$
It follows that $a_p\in\{-1,0,2\}$, and a theorem of Chebotarev
(generalising Dirichlet's theorem recalled above) implies that the proportion
of $p$ with $a_p=-1,0,2$ is $1/3, 1/2, 1/6$ respectively.

We thus have a ``formula'' for $N_p(T^3-T-1)$ similar in some sense to the one
for $N_p(T^2-T-1)$ above.  It is remarkable that $a_p={1\over2}(B_p-C_p)$
where $B_n$, $C_n$ are defined by the identities
$$
\sum_{j=0}^{+\infty}B_jq^j=\sum_{(x,y)\in\Z^2}q^{x^2+xy+6y^2},\quad
\sum_{j=0}^{+\infty}C_jq^j=\sum_{(x,y)\in\Z^2}q^{2x^2+xy+3y^2}.
$$
Even more remarkable is the fact that if we define a function of a complex
variable $\tau$ in ${\goth H}=\{x+iy\in\C\mid y>0\}$ by
$F(\tau)=\sum_{n=1}^{+\infty}a_ne^{2i\pi\tau.n}$, then $F$ is analytic and
``highly symmetric''~: for every matrix
% $$
% \pmatrix{a&b\cr c&d\cr}\quad (a,b,c,d\in\Z),\quad ad-bc=1,\quad
% c\equiv0\pmod{23},
% $$
$\pmatrix{a&b\cr c&d\cr}$ with $a,b,c,d\in\Z$, $ad-bc=1$, $c\equiv0\pmod{23}$, 
so that in particular the imaginary part of $\displaystyle {a\tau+b\over
  c\tau+d}$ is $>0$, we have 
$$
\displaystyle F\!\left({a\tau+b\over c\tau+d}\right)
=\left(d\over23\right)(c\tau+d)\,F(\tau),\quad
\left(d\over23\right)=
\cases{+1&if $\bar d\in\F_{23}^{\times2}$\phantom.\cr
 -1&if $\bar d\notin\F_{23}^{\times2}$.\cr}
$$
These symmetries of $F$ make it a ``primitive eigenform of weight~$1$,
level~$23$, and character $\left({}\over23\right)$''~: no mean achievement.
% Cf. Serre (2003) p.438, Koike (1984) p.87, 
% Ono (2004) p.18, Kioford (2008) p.101
This constitutes a ``reciprocity law'' for $T^3-T-1$.

\medbreak

There are ``reciprocity laws'' even for some polynomials $f\in\Z[S,T]$ in {\it
  two\/} indeterminates such as $f=S^2+S-T^3+T^2$~: there are nice
``formulae'' for the number $N_p(f)$ of zeros $(s,t)\in\F_p^2$ of $f$ (the
number of pairs $(s,t)\in\F_p^2$ satisfying $s^2+s-t^3+t^2=0$), and the
sequence of these numbers (or rather the numbers $a_p=p-N_p(f)$, which are in
the interval $[-2\sqrt p,+2\sqrt p]$ by a theorem of Hasse) is ``highly
symmetric'' as in the example of the polynomial $T^3-T-1$ above.  Let me
explain.

For our new $f=S^2+S-T^3+T^2$, defining $B_n$, $C_n$ by the identities
$$\eqalign{
\sum_{j=0}^{+\infty}B_jq^j
&=\sum_{(x,y,u,v)\in\Z^4}q^{x^2+xy+3y^2+u^2+uv+3v^2}\cr
&=1q^0+4q^1+4q^2+8q^3+20q^4+16q^5+32q^6+16q^7+\cdots\cr
&\;\;+4q^{11}+64q^{12}+40q^{13}+64q^{14}+56q^{15}+68q^{16}+40q^{17}+\cdots\cr
\sum_{j=0}^{+\infty}C_jq^j
&=\sum_{(x,y,u,v)\in\Z^4}q^{2(x^2+y^2+u^2+v^2)+2xu+xv+yu-2yv}\cr
&=1q^0+0q^1+12q^2+12q^3+12q^4+12q^5+24q^6+24q^7+\cdots\cr
&\;\;+0q^{11}+72q^{12}+24q^{13}+48q^{14}+60q^{15}+84q^{16}+48q^{17}+\cdots\cr
}
$$
we have $a_p={1\over4}(B_p-C_p)$ for every prime $p\neq11$.  Amazing~!  

This ``formula'' may look somewhat different from the previous one but it is
actually similar~: ${1\over4}(B_p-C_p)=c_p$ for every prime $p$ (and hence
$a_p=c_p$ for $p\neq11$), where $c_n$ is the coefficient of $q^n$ in the
formal product 
$$
\eta_{1^2,11^2}=q\prod_{k=1}^{+\infty}(1-q^k)^2(1-q^{11k})^2
=0+1.q^1+\sum_{n>1}c_nq^n.   
$$
We have already listed $c_p={1\over4}(B_p-C_p)$ for $p<19$~; here is
another sample~:
% $$
% \vbox{\halign{&\hfil$#$\quad\cr
% \multispan{15}\hrulefill\cr
% p=&2&3&5&7&11&13&17&19&23&29&31&37&41&\cdots\cr 
% \noalign{\vskip-5pt}
% \multispan{15}\hrulefill\cr
% c_p=&-2&-1&1&-2&1&4&-2&0&-1&0&7&3&-8&\cdots\cr
% \multispan{15}\hrulefill.\cr}}
% $$
$$
\vbox{\halign{&\hfil$#$\quad\cr
\multispan{12}\hrulefill\cr
p=&19&23&29&31&37&41&\cdots&1987&1993&1997&1999\cr 
\noalign{\vskip-5pt}
\multispan{12}\hrulefill\cr
c_p=&0&-1&0&7&3&-8&\cdots&-22&-66&-72&-20&\cr
\noalign{\vskip-5pt}
\multispan{12}\hrulefill.\cr}}
$$
As in the case of $T^3-T-1$, if we put $q=e^{2i\pi\tau}$ (with $\tau=x+iy$ and
$y>0$), we get an analytic function $F(\tau)=\sum_{n=1}^{+\infty}c_n
e^{2i\pi\tau.n}$  
% $$
% F(\tau)
% =e^{2i\pi\tau}\prod_{k=1}^{+\infty}
% (1-e^{2i\pi\tau.k})^2(1-e^{2i\pi\tau.11k})^2,
% $$
of $\tau\in{\goth H}$ which satisfies $\displaystyle F\!\left({a\tau+b\over
    c\tau+d}\right)=(c\tau+d)^2F(\tau)$ for every matrix $\displaystyle
\pmatrix{a&b\cr c&d\cr}$ ($a,b,c,d\in\Z$) with $ad-bc=1$ and
$c\equiv0\pmod{11}$.  In addition, $F$ has many other remarkable properties
(such as $c_{mm'}=c_mc_{m'}$ if $\gcd(m,m')=1$,
$c_{p^r}=c_{p^{r-1}}c_p-pc_{p^{r-2}}$ for primes $p\neq11$ and $r>1$, and
$c_{11^r}=c_{11}^r$ for $r>0$) giving it the elevated status of ``a primitive
eigenform of weight~$2$ and level~$11$''.

The analogy with the properties of the previous $F$ is perfect, with
$c\equiv0\pmod{11}$ replacing $c\equiv0\pmod{23}$ and $(c\tau+d)^2$ replacing
$({d\over23})(c\tau+d)$ in the expression for the value at
$(a\tau+b)/(c\tau+d)$.  (There is a sense in which the prime $11$ is ``bad''
for this $f=S^2+S-T^3+T^2$, much as $2$ is bad for $T^2+1$, $5$ is bad for
$T^2-T-1$, the prime divisors of $m$ are bad for $\Phi_m$, and $23$ is bad for
$T^3-T-1$.)

The above example is taken from Langlands and Harder, but perhaps goes back to
Eichler and implicitly to Hecke~; what we have said about $f$ (namely,
$a_p=c_p$ for all primes $p\neq11$) is summarised by saying that ``the
elliptic curve (of conductor~$11$) defined by $f=0$ is {\it modular\/}''.  It
can also be summarised by saying that the ``Artin $L$-function''
$$
L_1(s)=\prod_{p\neq11}{1\over 1-a_p.p^{s}+p.p^{-2s}}
$$
(with an appropriate factor for the prime~$11$) is the same as the ``Hecke
$L$-function'' $L_2(s)=\sum_{n>0}c_nn^{-s}$.  

There is a similar equality (up to factors at finitely many primes) of the
Artin $L$-function $\prod_p(1-a_p.p^{-s}+p.p^{-2s})^{-1}$ of the polynomial
$f=S^2-T^3+T$, defined using the number of zeros $N_p=p-a_p$ of $f$ modulo
various primes~$p$, with the Hecke $L$-function $\sum_{n>0}c_nn^{-s}$ of
$$
\eta_{4^2,8^2}=
q\prod_{k=1}^{+\infty}(1-q^{4k})^2(1-q^{8k})^2
=0+1.q^1+\sum_{n>1}c_nq^n
$$
the function $F(\tau)=\sum_{n>0}c_ne^{2i\pi\tau.n}$ corresponding to which
enjoys the symmetries
$$
\displaystyle F\!\left({a\tau+b\over c\tau+d}\right)
=(c\tau+d)^2F(\tau)
$$
for every $\tau\in{\goth H}$ (so $\tau=x+iy$ with $x,y\in\R$, $y>0$) and every
matrix 
$$
\pmatrix{a&b\cr c&d\cr}\quad (a,b,c,d\in\Z),\quad ad-bc=1,\quad
c\equiv0\pmod{32}.
$$
As a final example, consider $f=S^2-T^3-1$ with $N_p=p-a_p$ points
modulo~$p$ on the one hand, and 
$$
\eta_{6^4}=q\prod_{k>0}(1-q^{6k})^4=0+1.q^1+\sum_{n>1}c_nq^n
$$ 
on the other, for which the function $F(\tau)=\sum_{n>0}c_ne^{2i\pi\tau.n}$
has the above symmetries with $32$ replaced by $36$.  Each of these defines an
$L$-function, the Artin $L$-function $\prod_p(1-a_p.p^{-s}+p.p^{-2s})^{-1}$
and the Hecke $L$-function $\sum_{n>0}c_nn^{-s}$.  Remarkably, they are the
same (except for some adjustment at finitely many primes)~!

{\it Wiles and others (Taylor, Diamond, Conrad, Breuil) have proved that
  similar ``formulae'' and ``symmetries'' are valid --- such equalities of
  Artin $L$-functions with Hecke $L$-functions are valid --- for each of the
  infinitely many\/ $f\in\Z[S,T]$ which define an ``elliptic curve'' (such as,
  among many others tabulated by Cremona,
$$
f=S^2+S-T^3+T,\quad\hbox{\it and\/}\quad f=S^2+S -T^3-T^2+T+1,
$$
for which the primes\/ $37$ and $101$ are ``bad'' respectively), thereby
settling a conjecture of Shimura, Taniyama and Weil, and providing the
first proof of Fermat's Last Theorem in the bargain.}

The precise version of this uninformative statement goes under the slogan
``Every elliptic curve over $\Q$ is modular'' and has been called ``the
theorem of the century''~; it relies on the previous work of many
mathematicians in the second half of the XX$^{\rm th}$ century, among them
Langlands, Tunnell, Serre, Deligne, Fontaine, Mazur, Ribet, $\ldots$

\medbreak

\def\\{\oldstyle} 

The reciprocity law for abelian polynomials (class field theory,
$\\1900$--$\\1940$), the reciprocity law for the polynomial $T^3-T-1$ of group
${\goth S}_3$, the modulariy of the ``elliptic curve'' $S^2+S-T^3+T^2=0$ and
more generally the modulariy of every elliptic curve defined over $\Q$
($\\1995$--$\\2000$), and many other known reciprocity laws (Drinfeld,
Laumon-Rapoport-Stuhler, Lafforgue, Harris-Taylor, Henniart) which we have not
discussed, are all instances of the Langlands programme ($\\1967$--$+\infty?$).
At its most basic level, it is a search for patterns in the sequence $N_p(f)$
for varying primes $p$ and a fixed but arbitrary polynomial $f$ with
coefficients in $\Q$.  Our examples exhibit the ``pattern'' for some
specific~$f$.

Confining ourselves to monic irreducible polynomials $f\in\Z[T]$, Langlands
predicts that there is a ``reciprocity law'' for $f$ as soon as we give an
embedding $\rho:\Gal_f\rightarrow\GL_d(\C)$ of the group $\Gal_f$ into
the group $\GL_d(\C)$ of square matrices with complex entries and
determinant $\neq0$.  (What was the embedding for $f=T^3-T-1$~?  It came from
the unique irreducible representation of ${\goth S}_3$ on $\C^2$.)  Class
field theory is basically the case $d=1$ of this vast programme.  Essentially
the only other known instances of the precise conjecture are for $d=2$ and
either the group $\rho(\Gal_f)$ is solvable (Langlands-Tunnell) or the image
$\rho(c)$ of the element $c\in\Gal_f$ which sends $x+iy$ to $x-iy$ is not of
the form $\displaystyle\pmatrix{a&0\cr0&a}$ for any $a\in\C^\times$
(Khare-Wintenberger-Kisin).

\medskip

{\it So what is a reciprocity law~?\/} As the most basic level, it is the law
governing the primes modulo which a given monic irreducible $f\in\Z[T]$
factors completely, such as the quadratic reciprocity law when $f$ has
degree~$2$, or class field theory when $f$ is abelian.  We can also think of a
reciprocity law as telling us how many solutions a given system of polynomial
equations with coefficients in $\Q$ has over the various prime fields $\Fp$.
But we have seen that these laws can be expressed as an equality of an Artin
$L$-function with a Hecke $L$-function, and that the modularity of elliptic
curves over $\Q$ can also be expressed as a similar equality of $L$-functions.
We may thus take the view that a reciprocity law is encapsulated in the
equality of two $L$-functions initially defined in two quite different ways,
Artin $L$-functions coming from {\it galoisian representations\/} and Hecke
$L$-functions coming from {\it automorphic represetations}.

\bigskip
\centerline{---$*$---$*$---}
\bigbreak

The rest of this Note requires greater mathematical maturity.  We explain what
kind of ``reciprocity laws'' Serre's conjecture provides, but pass in silence
the recent proof of this conjecture by Khare and Wintenberger.

For any finite galoisian extension $K|\Q$ and any prime $p$ which does not
divide the discriminant of $K$, we have the ``element''
$\Frob_p\in\Gal(K|\Q)$, only determined up to conjugation in that group.  If
$f_p$ is the order of $\Frob_p$, then $p$ splits in $K$ as a product of
$g_p=[K_l:\Q]/f_p$ primes of $K$ of residual degree $f_p$, so that $p$ splits
completely if and only if $\Frob_p=\Id_K$.  The extension $K$ is uniquely
determined by the set $S(K)$ of primes which split completely in $K$.  For
example, if $S(K)$ is the set of $p\equiv1\pmod m$ for some given $m>0$, then
$K=\Q(\zeta_m)$, where $\zeta_m$ is a primitive $m$-th root of~$1$.  Given a
faithful representation of $\Gal(K|\Q)$ in $\GL_d(\C)$, Langlands'
programme characterises the subset $S(K)$ in terms of certain ``automorphic
representations of $\GL_d$ of the ad{\`e}les''.

Serre's conjecture deals with faithful representations $\rho$ of $\Gal(K|\Q)$
into $\GL_2(\bar\F_l)$ which are odd in the sense that the image $\rho(c)$ of
the ``complex conjugation'' $c\in\Gal(K|\Q)$ has determinant $-1$.  The oddity
is automatic if $l=2$~; for odd $l$, it implies that $K$ is totally imaginary.

A rich source of such representations is provided by elliptic curves $E$ over
$\Q$.  For every prime $l$, we have the $2$-dimensional $\F_l$-space
$E[l]={}_lE(\bar{\Q})$ of $l$-torsion points of $E$ together with a continuous
linear action of the profinite group $G_{\Q}=\Gal(\bar{\Q}|\Q)$.  Adjoining
the $l$-torsion of $E$ to $\Q$ we get a finite galoisian extension $\Q(E[l])$
which is unramified at every prime $p$ prime to $Nl$, where $N$ is the
conductor of $E$, and comes with an inclusion of groups
$\Gal(\Q(E[l])|\Q)\subset\Aut_{\F_l}({}_lE(\bar\Q))$ (which is often an
equality).  Choosing an $\F_l$-isomorphism $E[l]\to\F_l^2$, we
obtain a faithful representation $\rho_{E,l}:\Gal(\Q(E[l])|\Q)\to\GL_2(\F_l)$
which is odd because $\det\rho_{E,l}$ is the mod-$l$ cyclotomic character.

Given a finite galoisian extension $K$ of $\Q$ and a faithful representation
$\rho:\Gal(K|\Q)\to\GL_2(\bar\F_l)$ which is odd and irreducible (for the sake
of simplicity), Serre defines two numbers $k_\rho$, $N_\rho$ and conjectures
that $\rho$ ``comes from a primitive eigenform of weight $k_\rho$ and level
$N_\rho$''.  Instead of trying to explain this, let me illustrate it with the
example of the elliptic curve $E:y^2+y=x^3-x^2$ which we have already
encountered.

This curve has good reduction at every prime $p\neq11$.  Serre shows that the
representation $\rho_{E,l}:\Gal(\Q(E[l])|\Q)\to\GL_2(\F_l)$ is an isomorphism
for every $l\neq5$.  For every $p\neq11,l$, the characteristic polynomial of
$\rho_{E,l}(\Frob_p)\in\GL_2(\F_l)$ is
$$
X^2-\bar c_pX+\bar p\in\F_l[X]
$$
where $c_n$ is the coefficient of $q^n$ in
$\eta_{1^2,11^2}=q\prod_{k>0}(1-q^k)^2(1-q^{11k})^2$ as above.
Recall that $p$ splits completely in $\Q(E[l])$ if and only if
$$
\rho_{E,l}(\Frob_p)=\pmatrix{1&0\cr0&1},
$$
which implies, but is not equivalent to, $p,c_p\equiv1,2\pmod l$.  For
$l=7$, the first ten $p$ satisfying these congurences are 
$$
113, 379, 701, 1051, 2437, 2521, 2731, 2857, 3221, 3613~;
$$ 
none of them splits completely in $\Q(E[7])$.  Nevertheless, the
representation $\rho_{E,l}$ is explicit enough to determine the splitting of
rational primes in $\Q(E[l])$~; the first ten such primes for $l=7$ are
$$
4831, 22051, 78583, 125441, 129641, 147617, 153287, 173573, 195581,
$$ 
and $199501$, as obligingly computed by Tim Dokchitser at my request.

In summary, we have the following ``reciprocity law" for $\Q(E[l])$~:
$$
\hbox{``\thinspace $p\neq11,l$ splits completely in $\Q(E[l])$\thinspace''}
\quad\Leftrightarrow\quad
\hbox{``\thinspace}E_p[l]\subset E_p(\F_p)\hbox{\thinspace''},
$$
where $E_p$ is the reduction of $E$ modulo $p$.  Indeed, reduction modulo
$p$ identifies $E[l]$ with $E_p[l]$ and the action of ${\rm Frob}_p$ on the
former $\F_l$-space with the action of the canonical generator
$\varphi_p\in{\rm Gal}(\bar\F_p|\F_p)$ on the latter $\F_l$-space.  To say
that $\varphi_p$ acts trivially on $E_p[l]$ is the same as saying that
$E_p[l]$ is contained in the $\F_p$-rational points of $E_p$.  The analogy
with the multiplicative group $\mu$ is perfect~:
$$
\hbox{``\thinspace$p\neq l$ splits completely in $\Q(\mu[l])$\thinspace''}
\quad\Leftrightarrow\quad
\hbox{``\thinspace}\mu_p[l]\subset \mu_p(\F_p)\hbox{\thinspace''}
$$
($\Leftrightarrow l\,|\,p-1\Leftrightarrow p\equiv1\pmod l$), where $\mu_p$ is
{\it not\/} the $p$-torsion of $\mu$ --- that would be $\mu[p]$ --- but the
reduction of $\mu$ modulo $p$, namely the multiplicative group over $\F_p$. 

With the proof by Khare and Wintenberger (2006--2009) of Serre's modularity
conjecture (1973--1987), we can write down such reciprocity laws for every
finite galoisian extension $K|\Q$ as soon as a faithful odd irreducible
representation $\rho:\Gal(K|\Q)\to\GL_2(\bar\F_l)$ is given.

\medbreak

{\bf Sources.}  The reciprocity laws for quadratic and cyclotomic polynomials
can be found in most textbooks on number theory~; they are also treated at
some length by Wyman \citer\wyman().  An advanced
introduction to class field theory is contained in the book {\it Algebraic
  Number Theory\/} \citer\ant() edited by Cassels and Fr{\"o}hlich (recently
reissued by the London Mathematical Society with a long list of corrections,
available on Buzzard's website).

The polynomial $T^3-2$ is discussed by Cox \citer\cox().  I learnt the example
$T^3-T-1$ from Emerton and later in a paper by Serre \citer\serre().  It is
also discussed in an earlier paper \citer\serremodi(), as is the polynomial
$T^3+T-1$, for which the corresponding modular form is
$$
{1\over2}\left(\sum_{(x,y)\in\Z^2} q^{x^2+xy+8y^2}
   -\sum_{(x,y)\in\Z^2} q^{2x^2+xy+4y^2}\right).
$$
He also mentions a different kind of example, given by the dihedral extension 
$\Q(\root4\of1,\root4\of{12})$ of $\Q$, for which the corresponding modular
form is
$$
\eta_{12^2}=q\prod_{k=1}^{+\infty}(1-q^{12k})^2
 =\sum_{x,y\equiv1,0(3), x+y\equiv0(2)}(-1)^yq^{x^2+y^2}.
$$

The elliptic curve $E:y^2+y=x^3-x^2$ and its associated modular form
$\eta_{1^2,11^2}$ are mentioned by Langlands ({\it A little bit of number
  theory}, unpublished but available on the website of the Institute for
Advanced Study), and discussed in his talk \citer\gibbs() at the Gibbs
symposium (also available on the same website).  It is also discussed by
Harder \citer\harder() and available on the website of the International
Centre for Theoretical Physics), and by Zagier \citer\zagier().

This is the curve (with a different defining equation) used by Shimura
\citer\shimura() to study reciprocity laws for $\Q(E[l])$ for primes
$l\in[7,97]$.  Serre showed (Inventiones {\bf 15}) that $\rho_{E,l}$ is
surjective for every $l\neq5$~; cf.~the online notes on Serre's conjecture by
Ribet and Stein.  The other two elliptic curves, namely $y^2=x^3-x$ and
$y^2=x^3+1$, of conductor $32$ and $36$ respectively, are taken from a paper
by Honda \citer\honda().

\medbreak

{\bf Resources.}  In addition to the sources cited above, a number of articles
providing an introduction to the Langlands Programme are freely available
online.  There is a very accessible account by Taylor in \citer\taylor().
There are introductory accounts by Gelbart \citer\gelbart(), Arthur
\citer\arthur(), \citer\arthurmmiii(), Knapp \citer\knapp(), and the more
informal essays by Arthur, Knapp and Rogawski in the {\it Notices\/} of the
American Mathematical Society.  Knapp also has a number of other related
introductory texts on his website.  For the more advanced reader, there is a
book by Bump, Cogdell, de Shalit, Gaitsgory, Kowalski, and Kudla {\it An
  introduction to the Langlands program}, edited by Bernstein and Gelbart,
Birkhäuser (2003).

\medbreak

{\bf Applications.}  Reciprocity laws are a source of intense pleasure for the
mathematician~; they are also immensely useful within mathematics.  For the
rest, Simon Singh says that ``There are also important implications for the
applied sciences and engineering.  Whether it is modelling of the interactions
between colliding quarks or discovering the most efficient way to organise a
telecommunications network, often the key to the problem is performing a
mathemtical calculation.  In some areas of science and technology the
complexity of the calculations is so immense that progress in the subject is
severly hindered.  If only mathematicians could prove the linking conjectures
of the Langlands programme, then there would be shortcuts to solving
real-world problems, as well as abstract ones'' ({\it Fermat's Last Theorem},
p.~214). 

\medbreak

{\bf A final example.}  Let me end with another striking application of the
Langlands programme~: it allows us to prove some arithmetical statements which
have a fairly elementary formulation but for which no other proof is
available, elementary or otherwise.  Let $c_n$ ($n>0$) be the coefficient of
$q^n$ in the formal product
$$
\eta_{1^{24}}=
q\prod_{k=1}^{+\infty}(1-q^{k})^{24}=0+1.q^1+\sum_{n>1}c_nq^n.
$$
In $\oldstyle 1916$, Ramanujan had made some deep conjectures about these
$c_n$~; some of them (such as $c_{mm'}=c_mc_{m'}$ if $\gcd(m,m')=1$ and
$$
c_{p^r}=c_{p^{r-1}}c_p-p^{11}c_{p^{r-2}}
$$ 
for $r>1$ and primes~$p$) were proved by Mordell in $\oldstyle 1917$ but the
last of which was proved by Deligne only in the $\oldstyle70$s~: for every
prime $p$, the number $t_p=c_p/2p^{11/2}$ lies in the interval $[-1,+1]$.

All these properties of the $c_n$ follow from the fact that the corresponding
function $F(\tau)=\sum_{n>0}c_ne^{2i\pi\tau.n}$ of $\tau=x+iy$ ($y>0$) in
$\goth H$ is a ``primitive eigenform of weight~$12$ and level~$1$'' (which
basically amounts to the identity $F(-1/\tau)=\tau^{12}F(\tau)$).  Here are
the first few $c_p$ (usually denoted $\tau(p)$, but our $\tau$ is in $\goth
H$) computed by Ramanujan~:
$$
\vbox{\halign{&\hfil$#$\quad\cr
\multispan{9}\hrulefill\cr
p=&2&3&5&7&11&13&17&\cdots\cr
\noalign{\vskip-5pt}
\multispan{9}\hrulefill\cr
c_p=&-24&252&4830&-16744&534612&-577738&-6905934&\cdots\cr
\noalign{\vskip-5pt}
\multispan{9}\hrulefill.\cr}}
$$
(Incidentally, Ramanujan had also conjectured some congruences satisfied by
the $c_p$ modulo $2^{11}$, $3^7$, $5^3$, $7$, $23$ and $691$, such as
$c_p\equiv1+p^{11}\pmod{691}$ for every prime $p$~; they were at the origin of
Serre's conjecture mentioned above.)

We may therefore ask how these $t_p=c_p/2p^{11/2}$ are distributed~: for
example are there as many primes $p$ with $t_p\in[-1,0]$ as with
$t_p\in[0,+1]$~?  Sato and Tate predicted in the $\\60$s that the precise
proportion of primes $p$ for which $t_p\in[a,b]$, for given $a<b$ in
$[-1,+1]$, is
$$
{2\over\pi}\!\int_{\!a}^b\!\!\!\sqrt{1-x^2}\;dx.
$$
This is expressed by saying that the $t_p=c_p/2p^{11/2}$ are {\it
  equidistributed\/} in the interval $[-1,+1]$ with respect to the measure
$(2/\pi)\sqrt{1-x^2}\;dx$.  Recently Barnet-Lamb, Geraghty, Harris and Taylor
have proved that such is indeed the case.

Their main theorem implies many such equidistribution results~; for example,
going back to the equation $S^2+S=T^3-T^2$, which has $p-a_p$ solutions in
$\F_p^2$, where $a_p$ is the coefficient of $q^p$ in $\eta_{1^2,11^2}$ for
every prime $p\neq11$, the numbers $a_p/2\sqrt p$, which lie in $[-1,+1]$,
are equidistributed in this interval with respect to the same measure.  See
\citer\taylor() for an introduction to such {\it density theorems}.

\medbreak

{\bf Acknowledgments.}  I heartily thank Tim Dokchitser for computing the
first ten primes which split completely in $\Q(E[7])$, and my colleague Rajesh
Gopakumar for a careful reading of the text.

\bigbreak
\unvbox\bibbox 
\bigbreak
%\vfill\eject
\vskip2cm plus1cm minus1cm
{\obeylines\parskip=0pt\parindent=0pt
Chandan Singh Dalawat
Harish-Chandra Research Institute
Chhatnag Road, Jhunsi
{\pc ALLAHABAD} 211\thinspace019, India
\vskip5pt
\tt dalawat@gmail.com}

\bye